\documentclass[11pt]{article} 

\usepackage{fullpage}
\usepackage{natbib}
\usepackage{amsthm}
\usepackage{color}
\usepackage{amsmath}
\usepackage{amsfonts}
\usepackage{xcolor}
\colorlet{linkequation}{blue}
\usepackage{framed}
\usepackage{times}
\usepackage{amsmath}
\usepackage{amsthm}
\usepackage{amssymb}
\usepackage{wasysym}
\usepackage{mathrsfs}
\usepackage{bbm}
\usepackage{geometry} 
\usepackage{mathtools}
\usepackage[colorlinks=true,citecolor=blue,linkcolor=blue]{hyperref}

\theoremstyle{plain}
\newtheorem{theorem}{Theorem}
\newtheorem{lemma}{Lemma}
\newtheorem{definition}{Definition}

\newtheorem{corollary}{Corollary}
\newtheorem{remark}{Remark}

\newcommand{\E}{{\mathbb E}}

\newcommand{\R}{{\mathbb R}}
\renewcommand{\P}{{\mathbb P}}

\newcommand{\C}{{\mathcal{C}(1)}}

\newcommand{\F}{{\cal F}}

\newcommand{\PP}{\mathbb{P}}

\newcommand{\eps}{\epsilon}
\def\qt#1{\qquad\text{#1}}

\newcommand{\Unif}{\P}

\newcommand{\Sd}{\mathbb{S}^{d-1}}
\newcommand{\Ri}{\mathcal{R}}

\newcounter{rcnt}[section]

\def\qt#1{\qquad\text{#1}}
\newcommand{\hp}{h_{P_{m,d}}}
\newcommand{\V}{\mathcal{V}}
\def\argmin{\mathop{\rm argmin}}
\def\argmax{\mathop{\rm argmax}}
\newcommand{\nnt}{n^{-2/(d-1)}}
\newcommand{\mmt}{m^{-2/(d-1)}}

\newcommand{\GW}{\bar{\mathcal{W}}}

\newcommand{\Fcal}{\mathcal{F}}
\newcommand{\Ncal}{\mathcal{N}}

\newcommand{\cG}{{\mathcal G}}

\usepackage{times}

\begin{document}
\title{On Suboptimality of Least Squares with Application to Estimation of Convex Bodies}

\author{Gil Kur \thanks{Massachusetts Institute of Technology, \texttt{gilkur@mit.edu}}
	\and
	Alexander Rakhlin \thanks{Massachusetts Institute of Technology, \texttt{rakhlin@mit.edu}}
	\and
	Adityanand  Guntuboyina \thanks{UC Berkeley, \texttt{aditya@stat.berkeley.edu}.}}
\date{}
\maketitle

\begin{abstract}%
	We develop a technique for establishing lower bounds on the sample complexity of Least Squares (or, Empirical Risk Minimization) for large classes of functions. As an application, we settle an open problem regarding optimality of Least Squares in estimating a convex set from noisy support function measurements in dimension $d\geq 6$. Specifically, we establish that Least Squares is mimimax \textit{sub-optimal}, and achieves a rate of $\tilde{\Theta}_d(n^{-2/(d-1)})$ whereas the minimax rate is $\Theta_d(n^{-4/(d+3)})$. 
\end{abstract}

\paragraph{keywords}%
	Non-parametric statistics, regression, support function, ERM, Least squares, non-Donsker regime

\section{Introduction and main results}
Consider the problem of regression where the goal is to estimate a function $f^* : \mathcal{X} \rightarrow \R $  from observations $(X_1, Y_1), \ldots, (X_n, Y_n)$ drawn according to the model 
\begin{equation*}
Y_i = f^*(X_i) + \xi_i \qt{for $i = 1, \dots, n$}
\end{equation*}
where $X_1, \dots, X_n$ are design points (fixed or random) and $\xi_1, \dots, \xi_n$ are i.i.d random variables with mean zero.  Assume that the unknown function $f^*$ belongs to a known \emph{convex} class of functions $\mathcal{F}$. 
The most natural and basic estimator in this setting is the Least Squares Estimator (LSE)
\begin{equation*}
\hat{f}_n \in \argmin_{f \in \mathcal{F}} \sum_{i=1}^n \left(Y_i - f(X_i) \right)^2, 
\end{equation*}
which is also the Maximum Likelihood Estimator (MLE) under Gaussian noise. It is customary to evaluate the performance of $\hat{f}_n$ via the risk 
\begin{equation}\label{fdloss}
\|\hat{f}_n - f^*\|^2_{\PP_n} :=  \frac{1}{n} \sum_{i=1}^n \left(\hat{f}_n(X_i) - f^*(X_i) \right)^2
\end{equation}
in the fixed-design setting where $X_1, \dots, X_n$ are assumed to be non-random points in $\mathcal{X}$, and via 
\begin{equation}\label{rdloss}
\|\hat{f}_n - f^*\|^2_{\PP} := \int_{\mathcal{X}} \left(\hat{f}_n(x) - f^*(x) \right)^2 d\PP(x)
\end{equation}
in the random design setting where $X_1, \ldots, X_n \overset{i.i.d}{\sim} \PP$. We shall focus on the optimality of the LSE as measured in the minimax sense (see e.g., \cite{tysbakovnon}) via the worst-case risk 
\[
\Ri_{\PP_n}(\hat{f}_n,\F) := \sup_{f^* \in \F}\E_{f^*} \|\hat{f}_n - f^*\|_{\PP_n}^2 ~~~ \text{ or } ~~~ \Ri_{\PP}(\hat{f}_n,\F) := \sup_{f^* \in \F}\E_{f^*} \|\hat{f}_n - f^*\|_{\PP}^2,
\] 
in the fixed and random designs, respectively. 

The accuracy of the LSE is one of the most basic questions in statistics \citep{van2000empirical,birge1993rates,koltchinskii2011oracle,van1996weak}. Typical results impose conditions on the metric entropy on (natural subsets of) the class $\mathcal{F}$. We shall denote the $\epsilon$-metric entropy of a class $\mathcal{G}$ of functions under the $L^2(\mathbb{Q})$ pseudometric by $\log \Ncal(\epsilon, \mathcal{G}, \mathbb{Q})$. 
Here $\Ncal(\epsilon, \mathcal{G}, \mathbb{Q})$ is the minimal cardinality of a set $\mathcal{N}_{\epsilon}$ of functions such that for any $g\in\mathcal{G}$ there exists a $g'\in\mathcal{N}_{\epsilon}$ with $\|g-g'\|_{\mathbb{Q}} \leq \epsilon$.

It is convenient to isolate the existing results on the LSE into two categories depending on whether the class of functions $\mathcal{F}$ is in the \textit{Donsker} regime or in the \textit{non-Donsker} regime. By the Donsker (or, more precisely, $\mathbb{P}$-Donsker) regime, we mean that the $L^2(\mathbb{P})$ metric entropy of natural totally bounded subsets of $\mathcal{F}$ grows as $\epsilon^{-p}$ for some $p < 2$.\footnote{In this section, for brevity we may omit constants that may depend on $p$.} In contrast, for the non-Donsker regime, the entropy grows as $\epsilon^{-p}$ for some $p > 2$. 

It should be noted that the two regimes can often occur in different settings of the same general problem. For example, the class of support functions of compact, convex subsets of $\R^d$ is Donsker for $d \leq 4$ and non-Donsker for $d \geq 6$. We shall revisit this example in detail in Section \ref{suppest}. 

In the Donsker regime, it is well-known that the LSE achieves the rate $n^{-2/(2+p)}$, which is the minimax rate of convergence, under minimal assumptions \citep{barron1999risk}. Thus, the LSE is minimax optimal when the class $\mathcal{F}$ is Donsker with respect to $\mathbb{P}$ or $\mathbb{P}_n$. In contrast, the minimax optimality of the LSE in the non-Donsker regime remains unresolved. In a fundamental paper, \cite{birge1993rates} proved that in the non-Donsker regime, the rate of convergence of the LSE is bounded from above by $n^{-1/p}$. However, it is also well-known (see e.g., \cite{yang1999information}) that the minimax rate of estimation is still $n^{-2/(2+p)}$. 

It was first observed in \citep[Section 4]{birge1993rates} and later in \cite{birge2006model} that it is possible to design ``pathological'' function classes $\mathcal{F}$ where the LSE provably achieves a risk of $\Tilde{\Theta}(n^{-1/p})$. However, for a general non-Donsker class $\mathcal{F}$, the fundamental question of whether the rate of convergence of the LSE is $n^{-1/p}$ or $n^{-2/(2+p)}$ (or some intermediate rate) is unresolved.

Some important progress on this open problem has been made in the recent papers \citep{kur2019optimality,han2019global,carpenter2018near,han2019isotonic}. Specifically, these papers have shown that there exist ``natural'' non-Donsker families of functions (including the class of bounded convex functions on a smooth domain in $\R^d$ for $d \geq 3$, the class of multivariate Isotonic functions over $\R^d$ for $d \geq 2$, etc.) where the LSE achieves the minimax rate of $n^{-2/(2+p)}$, in contrast to the pessimistic upper bound of $n^{-1/p}$. 

In this paper, we give a recipe for establishing a lower bound for the LSE's risk for convex families of functions in the \emph{non}-Donsker regime both in fixed and random design. As an application, we complement the aforementioned recent results by proving that there also exist other ``natural'' non-Donsker families of functions (such as the class of support functions of compact convex subsets of $\R^d$ for $d \geq 6$) where the LSE cannot achieve a rate faster than $n^{-1/p}$, up to logarithmic multiplicative factors. In other words, for these non-Donsker classes, the LSE is provably \textit{suboptimal}. 

We start by stating assumptions that imply sub-optimality of LSE in a fixed design. Then, we extend these assumption to the random design setting. These assumptions will be verified for the case of support functions in Section~\ref{suppest}. 

Throughout, we assume that $\mathcal{F}$ is a \emph{convex} family of functions. 
Also, in the remainder of the paper, we assume for simplicity that $\xi \sim N(0,1)$; see Remark \ref{Remark1} for more general assumptions on the noise. We employ the notation $a\lesssim b$ to mean inequality up to a constant that may depend on $p$. Let 
$$B_{\PP_n}(f_0, t) := \left\{f \in \mathcal{F} : \|f - f_0\|_{\PP_n} \leq t \right\}$$ and for a set of functions $\mathcal{G}$ on $\mathcal{X}$, the Gaussian width of $\mathcal{G}$ is defined, conditionally on $X_1,\ldots,X_n$, as
\begin{equation}\label{Eq:GWdef}
\GW(\mathcal{G}) := \E \sup_{f \in \mathcal{G}} \frac{1}{n}
\sum_{i=1}^n g_i  f(X_i),
\end{equation}
where $g_1,\ldots,g_n$ are i.i.d. $N(0,1)$. 

We are now ready to state our key structural assumptions.
\paragraph{Assumptions: Fixed design} Let $p > 2$  and $n \gtrsim 1$. Assume there are two functions $f_0,f_{(n)} \in \Fcal$ such that for some $\beta, \gamma \geq 0$,
\begin{equation}\label{A21}
\log \Ncal(\epsilon, B_{\PP_n}(f_0,2\epsilon), \PP_n) \gtrsim \epsilon^{-p} \quad   \forall~ \epsilon  \gtrsim n^{-1/p}\log(n)^{\gamma}.
\end{equation}
Furthermore, suppose
\begin{equation}\label{A31}
\|f_{(n)} - f_0\|_{\PP_n} \lesssim n^{-\frac{1}{2p}}\log(n)^{-s_{p,\gamma}},
\end{equation}
where $s_{p,\gamma} = \frac{\gamma}{2}(\frac{p}{2}-1)$, and either 
\begin{equation}\label{A41}
\log \Ncal(\epsilon, B_{\PP_n}(f_{(n)}, t), \PP_n) \lesssim \sqrt{n}\log(n)^{p \cdot s_{p,\gamma}}\left(\log \frac{1}{\epsilon} \right)^{\beta} \left(\frac{t}{\epsilon} \right)^p
\end{equation}
for every $ n^{-\frac{1}{p}}\log (n)^{-s_{p,\gamma}} \lesssim \eps \leq t \lesssim n^{-\frac{1}{2p}} \log (n)^{-\frac{\beta}{p} -\frac{3s_{p,\gamma}}{2}}$ or a weaker condition
\begin{align}
\label{eq:weaker_cond}
\GW(B_{\PP_n}(f_{(n)}, t)) \lesssim n^{-\frac{1}{2p}}\log(n)^{\frac{\beta}{p} + s_{\gamma,p}}t,
\end{align}
for some $t \lesssim n^{-\frac{1}{2p}} \log (n)^{-\frac{\beta}{p} -\frac{3s_{p,\gamma}}{2}}$ holds.

The following theorem establishes that the rate of convergence of the LSE equals $n^{-1/p}$, up to a poly-logarithmic multiplicative factor, under the aforementioned assumptions. 
\begin{theorem}\label{Thm:LowerBirgeMassartFixed}
	Assume that $\mathcal{F}$ is a convex class of functions satisfying \eqref{A21}, \eqref{A31}, and either \eqref{A41} or \eqref{eq:weaker_cond} holds for some $p > 2$ and $n \gtrsim 1$. Then the following holds:
	\begin{equation*}
	\mathcal{R}_{\PP_n}(\hat{f}_n, \F) \gtrsim n^{-\frac 1p}\log(n)^{-\frac{2\beta}{p} -3s_{p,\gamma}}.
	\end{equation*}
\end{theorem}
We shall next extend Theorem \ref{Thm:LowerBirgeMassartFixed} to the case of random design, where the design points $X_1, \dots, X_n$ are independent with common distribution $\PP$. It is sufficient for our purposes to establish near-isometry: with high probability,
\begin{equation}\label{Eq:Distortion}
\begin{aligned}
\forall f,g\in\F,~~~ \frac{1}{2}\|f-g\|_{\PP_n} - d_n \leq \|f-g\|_{\PP} \leq  2\|f-g\|_{\PP_n} + d_n
\end{aligned}
\end{equation}
for some remainder $d_n$ that decays to zero with increasing $n$. In the case of uniformly bounded functions, sufficient conditions for the two-sided inequality \eqref{Eq:Distortion} can be found in the literature on local Rademacher averages (see e.g., \citep{bousquet2002concentration,bartlett2005local}), while the right-hand side of \eqref{Eq:Distortion} holds under weaker conditions \citep{koltchinskii2013bounding,mendelson2014learning,mendelson2017aggregation}. For the purpose of proving lower bounds for random design, however, we need the more demanding left-hand side of \eqref{Eq:Distortion}. Hence, we shall assume that functions in $\Fcal$ are uniformly bounded.

In addition, we assume the growth of Koltchinskii-Pollard entropy 
\begin{equation}\label{A1}
\sup_{n \in \mathbb{N}}\sup_{{\mathbb Q} \in \mathcal{P}_n}\log \mathcal{N}(\epsilon,\F,{\mathbb Q})   \lesssim \epsilon^{-p} 
\end{equation}
where $\mathcal{P}_{n}$ denotes the set of all probability measures supported on finite subsets of $\mathcal{X}$ of cardinality at most $n$. Under this assumption, with high probability $d_n \lesssim(\log n)^{2} n^{-1/p}$  \citep{rakhlin2017empirical}. This will allow us to reduce the random design setting to a fixed design and use Theorem \ref{Thm:LowerBirgeMassartFixed}, since the remainder $(d_n)^2$ is of the lower order than $n^{-1/p}$.\footnote{For some families of functions it might be possible to relax the additional assumptions of uniform boundedness and \eqref{A1} and achieve a similar result to Corollary \ref{Thm:LowerBirgeMassartRandom}, see for example \cite{KurConv}} 

To establish a lower bound for LSE in random design, we may verify that above assumptions \eqref{A21}, \eqref{A31}, \eqref{A41} hold with high probability for random measures $\PP_n$ and employ near-isometry for the distance between LSE and the regression function. Alternatively, it may be easier to verify the corresponding assumptions in the population. We now state this latter approach. 

First, recall that bracketing number $\Ncal_{[]}(\epsilon, \mathcal{G}, \PP)$ is the minimal cardinality of a set $\mathcal{N}_{[],\epsilon}$ of pairs of functions $(g^{-},g^{+})$ such that $\|g^{+}-g^{-}\|_{\PP} \leq \epsilon$, and for any $g \in \mathcal{G}$ there exists $(g^{-},g^{+}) \in  \mathcal{N}_{[],\epsilon}$ such that $g^{-} \leq g \leq g^{+}.$ 

\paragraph{Assumptions: Random design} Let $p > 2$. Assume there exists a function $f_0 \in \Fcal$ such that for every integer $m \gtrsim 1$, there exists $f_{(m)} \in \mathcal{F}$ 
such that
\begin{equation}\label{A2}
\log \Ncal(\epsilon, B_{\PP}(f_0,2\epsilon), \PP) \gtrsim \epsilon^{-p} \quad   \forall \epsilon  \geq  0
\end{equation} 
and
\begin{equation}\label{A3}
\|f_{(m)} - f_0\|_{\PP} \lesssim m^{-\frac{1}{p}}.
\end{equation}
Furthermore, for some $\beta \geq 0$
\begin{equation}\label{A4}
\log \Ncal_{[]}(\epsilon, B_{\PP}(f_{(m)}, t), \PP) \lesssim m \left(\log \frac{1}{\epsilon} \right)^{\beta} \left(\frac{t}{\epsilon} \right)^p,
\end{equation}
for every $ 0< \epsilon \leq t \lesssim \log(m)^{-\beta/p}m^{-1/p}$. 

We remark that assumption \eqref{A2} is satisfied for non-Donsker classes of functions: see \cite[Lemma 3]{yang1999information} for the proof of existence. 

The following is our lower bound for the random-design setting:
\begin{corollary}\label{Thm:LowerBirgeMassartRandom}
	Assume that $\F$ is a convex class of uniformly bounded functions satisfying \eqref{A1}, \eqref{A2}, \eqref{A3}, and \eqref{A4} for some $p>2$ and $n\gtrsim 1$. Then
	\begin{equation*}
	\mathcal{R}_{\PP}(\hat{f}_n, \F) \gtrsim n^{-\frac 1p}\log(n)^{-\frac{2\beta}{p} -3(\frac{p}{2}-1)}.
	\end{equation*}
\end{corollary}  

\subsection{The ideas behind our assumptions}

Assumption \eqref{A2} says that the function $f_0$ is ``complex'' in that its local neighborhood is as complex as the entire function class $\mathcal{F}$ when $\log\Ncal(\epsilon,\F,\PP) \sim \epsilon^{-p}$. 
The assumption \eqref{A3} states that $f_0$ is approximated by $f_{(m)}$ up to the accuracy $m^{-1/p}$ for every $m$. The assumption \eqref{A4} captures the ``simplicity'' of the functions $f_{(m)}$ in relation to the complex function $f_0$ satisfying \eqref{A2}. Note that when $t = O(\epsilon)$, the right hand side of \eqref{A4} is logarithmic in $\epsilon$ (assuming that $m$ is not too large) while the right hand side of \eqref{A2} is polynomial in $1/\epsilon$. Thus the local neighborhood of $f_{(m)}$ is smaller than that of $f_0$ and in this sense $f_{(m)}$ is simpler than $f_0$. Note also that there is a factor of $m$ on the right hand side of \eqref{A4} which means that the complexity of the functions $f_{(m)}$ increases with $m$.

There exist natural non-Donsker function classes $\mathcal{F}$ which satisfy \eqref{A2}, \eqref{A3} and \eqref{A4}. In Section \ref{suppest}, we show that the class of support functions of compact convex sets in $\R^d$ satisfy these assumptions with $f_0$ being the support function of the unit ball and $f_{(m)}$ the support function of a regular polytopal approximation to the unit ball with $m$ vertices. 

The proof of Theorem \ref{Thm:LowerBirgeMassartFixed} will reveal that $\hat{f}_n$ achieves the rate $n^{-1/p}$ (up to logarithmic factors) when the regression function is in the class $\{f_{(m)}, m \geq 1\}$. In fact, the specific function achieving the rate $n^{-1/p}$ is $f^*=f_{(m)}$ for $m \sim \sqrt{n}$. It is interesting to note that for $m \sim \sqrt{n}$ and $t^2 \sim n^{-1/p}$, the local neighborhood $B_{\PP_n}(f_{(m)}, t)$ has the same metric entropy as the right hand side of \eqref{A2}. Our main technical insight is that the sub-optimality occurs at a function $f^*=f_{(m)}$ instead of perhaps a more natural candidate function such as $f_0$ (we actually believe that the rate at $f_0$ may be equal to $n^{-2/(2+p)}$).

\subsection{Discussion: the non-Donsker regime, revisited}
From the recent results on families in the non-Donsker regime \citep{kur2019optimality,han2019global,carpenter2018near,KurConv,han2019isotonic} and the present results, the LSE may be optimal or sub-optimal in the non-Donsker regime for natural classes of functions, even if uniformly bounded. These results also indicate that LSE achieves a risk that equals to the Gaussian complexity of the class, up to a constant that depends on dimension. Namely, in contrast to the Donsker regime, there is no localization. In the problem of convex regression with uniformly bounded functions and  \emph{Euclidean ball} as domain, as well as in multiple isotopic regression, the minimax rate equals to the Gaussian complexity of the family. In contrast, by Theorem \ref{Cor:SuppRandom}, for support function regression and for convex uniformly bounded regression (or Lipshitz-convex regression) with support on the \emph{cube} \citep{KurConv}, the Guassian complexity differs from the minimax rate.

\section{Application: Sub-optimality of Least Squares for estimating a convex set from noisy support function measurements}\label{suppest}
In this section, we shall use Corollary  \ref{Thm:LowerBirgeMassartRandom} to resolve a long-standing open problem on the sub-optimality of the least squares estimator in the problem of estimating a convex set in dimension $d\geq6$ from noisy support function measurements, see for example \citep{gardner2006convergence,brunel2016adaptive,brunel2013adaptive,brunel2018estimation,guntuboyina2012optimal,fisher1997estimation,balazs2015near,soh2019fitting}. Let us first recall that the support function $h_K: \Sd \rightarrow \R$ of a compact convex set $K$ in $\R^d$ is defined as
\begin{equation*}
h_K(u) = \max_{x \in K} \left<x, u \right> \qt{for $u \in \Sd := \left\{u : \|u\| = 1 \right\}$}. 
\end{equation*}
The support function uniquely determines the compact convex set $K$ and is a fundamental object in convex geometry (see, for example, \cite[Section 1.7]{schneider2014convex} or \cite[Section 13]{rockafellar1970convex}). Consider now the problem of estimating an unknown compact, convex set $K^*$ from observations $(X_1, Y_1), \dots, (X_n, Y_n)$ drawn according to the model: 
\begin{equation*}
Y_i = h_{K^*}(X_i) + \xi_i \qt{for $i = 1, \dots n$}
\end{equation*}
where $X_1, \dots, X_n$ are design points (fixed or random) and $\xi_1, \dots, \xi_n$ are i.i.d. $N(0,1)$. Recovery of $K^*$ is a fundamental problem in geometric tomography \citep{prince1990reconstructing,gardner1995geometric}. The natural estimator in this problem is the least squares estimator defined by 
\begin{equation*}
\hat{K}_n \in \argmin_{K \in \mathcal{C}} \sum_{i=1}^n \left(Y_i - h_K(X_i) \right)^2
\end{equation*}
where $\mathcal{C}$ denotes the class of all compact, convex sets in $\R^d$. Basic properties and algorithms for computing $\hat{K}_n$ can be found in \citep{prince1990reconstructing,lele1992convex} and \citep{kiderlen2008new}. Rigorous accuracy results for $\hat{K}_n$ as an estimator of $K^* \subset B_d$ (here $B_d$ denotes the unit-Euclidean ball) were proved in \citep{gardner2006convergence}. Specifically, \cite[Corollary 5.7]{gardner2006convergence} proved that, under the fixed design setting of any $X_1,\ldots,X_n$ (fixed points) that form a well-separated set, (i.e. a set of $n$ points on the unit-sphere that are at least $c(d)n^{-\frac{1}{d-1}}$ far from each other),
\begin{equation}\label{gardub}
\|h_{\hat{K}_n} - h_{K^*}\|_{\PP_n}^2 = O_P(\beta_n) \quad \text{ where } 
\beta_n = 
\begin{cases} 
n^{-4/(d+3)} & \text{for } d = 2, 3, 4 \\
n^{-4/(d+3)} \cdot \log n       & \text{for } d = 5 \\
n^{-2/(d-1)} & \text{ for } d \geq 6. 
\end{cases}
\end{equation}
Complementarily, \cite{guntuboyina2012optimal} proved that the minimax rate of estimation in this problem (under the same fixed design version of the problem and the $\|\cdot\|_{\PP_n}$ loss function) equals $\Theta_d(n^{-4/(d+3)})$ for all $d \geq 2$. These two results combined imply that the least squares estimator $\hat{K}_n$ is minimax optimal for $d = 2, 3, 4$ and nearly minimax optimal (up to the logarithmic multiplicative factor $\log n$) for $d = 5$. However, there is a gap between the upper bound $n^{-2/(d-1)}$ on the rate of convergence of the least squares estimator \eqref{gardub} and the minimax rate $n^{-4/(d+3)}$ for $d \geq 6$. This gap has remained open since the paper \cite{gardner2006convergence} where it was suggested that the upper bound is accurate and that the least squares estimator is indeed minimax \emph{suboptimal} for $d \geq 6$. The goal of this section is to confirm the conjecture of \cite{gardner2006convergence}. 

Specifically, we use Corollary  \ref{Thm:LowerBirgeMassartRandom} to prove that the LSE is suboptimal for $d \geq 6$, in the sense that there exist sets  for which the rate of convergence for LSE is bounded from below by $n^{-2/(d-1)}$ up to a logarithmic multiplicative factor. Moreover, in this result, $f_0$ will be the $d-$dimensional ball, and the sets $\{K_{m}\}_{m=d+1}^{\infty} \subset \mathcal{C}$ will be ``regular'' polytopes that form an optimal approximation to the ball with $m \geq d+1$ vertices (see Lemma \ref{Lem:Regu} for more details).

\begin{corollary}\label{Cor:SuppRandom}
	Let $d \geq 6$. Suppose $X_1, \dots, X_n \sim \PP,$ where $\PP$ is the uniform distribution\footnote{Or for any density $g(x)$ on the sphere, such that $g(x) \geq c_1$ for all $x \in \Sd$.} on $\Sd$. There exist a positive constants, $C_d$, $\gamma_d$ depending only on $d$ such that  
	\begin{equation*}
	\sup_{K^* \in \mathcal{C}(1)} \E_{K^*} \|h_{\hat{K}_n} - h_{K^*}\|^2_{\PP}  \geq C_d n^{-2/(d-1)}\log(n)^{-\gamma_d}, 
	\end{equation*}
	where $\mathcal{C}(1)$ is the set of all compact convex sets that are contained in the unit Euclidean ball in $\mathbb{R}^d$. 
\end{corollary}
Observe that here the LSE in not restricted to $\mathcal{C}(1)$, i.e. $h_{\hat{K}_n}$ can return any convex set in $\R^d$.

Finally, the proof of the above corollary can be modified to hold for any fixed design of $n$ well-separated points. The modification will follow from the fact that $n$ well-separated points are a discrete approximation to the uniform measure on the sphere. Therefore, we settle the question in \citep{gardner2006convergence}.

\paragraph{Acknowledgments}
The first and second authors acknowledge the support from NSF under award DMS-1953181, and from the Center for Minds, Brains and Machines (CBMM) funded by NSF award CCF-1231216. The third author is supported by NSF CAREER Grant DMS-1654589. We also acknowledge the anonymous referees for their useful comments.

\section{Proofs}

\paragraph{Notation}
Throughout this text, $c,C$ with subscripts are positive absolute constants that do not depend on the dimension $d$. Additionally, positive  constants  that only depend on the dimension or on $p$ are explicitly denoted, respectively, by $c(d),C(d),c(p), C(p)$. These constants may change from line to line, and from section to section. $B_d$ denotes the unit Euclidean ball in $\R^d$. $\Tilde{O}(\cdot),\Tilde{\Omega}(\cdot)$ denotes behavior up to logarithmic factors of ($n$ or $\epsilon$).

\subsection{Proof of Theorem \ref{Thm:LowerBirgeMassartFixed}}
Our main technical tool is the following important result of \cite{chatterjee2014new} which gives sharp upper and lower bounds for the accuracy of an LSE over a convex family of functions in the fixed-design setting. We use the following notation in this result.

\begin{theorem}[ Theorem 1.1 \cite{chatterjee2014new}]\label{chat}
	Let $\F$ be a convex family of functions and consider the LSE for the fixed design setting. Let 
	\begin{equation}\label{chat.t}
	t_f := \argmax_{t \geq 0}H_f(t)\qt{ where $H_f(t) := \GW(B_{\PP_n}(f,t)) - \frac{t^2}{2}$}.
	\end{equation}
	Then, $t_f$ is unique, $H_{f}(\cdot)$ is a concave function, and  we have 
	\begin{equation}\label{Eq:SouarvBd}
	\P \left\{0.5t_f^2\leq \|\hat{f}_n - f\|^2_{\PP_n}  \leq 2t_f^2 \right\} \geq 1 - 3p_n,
	\end{equation}
	where $p_n = \exp \left(-c n t_f^2 \right)$.
\end{theorem}
\begin{remark}\label{Remark1}
	Based on the proof of \cite[Thm 1.1]{chatterjee2014new}, the assumption on $g = (g_1,\ldots,g_n) \sim N(0,I_{n \times n})$ can be relaxed to any isotropic random vector that satisfies a convex Lipshitz concentration inequality \citep{adamczak2015note,boucheron2013concentration}.  Also, the same result for an isotropic log-concave noise with $p_n= \exp \left(-c \sqrt{n} t_f^2 \right)$ follows from the recent result of \cite{lee2018stochastic}, which is a concentration inequality for Lipshitz functions for isotropic log-concave distributions. Finally, we remark that for zero-mean $\xi$ with variance $\sigma^2$, all our results scale by $\sigma^2$.
\end{remark}
The proof of Theorem \ref{Thm:LowerBirgeMassartFixed} is reduced to the following ``two points'' lemma that is mainly based on Theorem \ref{chat}. The notation $f_{(n)}$ is to emphasize the fact that $f_{(n)}$ is a function that is chosen based on the number of samples.
\begin{lemma}\label{Lem:Twoptsfixed}[Lower bound for fixed design]
	Let $f_0,f_{(n)} \in \F$, and $r_n \geq Cn^{-1/4},\delta_n,w_n,s_n$ be positive constants. Also, let $w_n' = \max\{\frac{r_n^2}{8\delta_n},w_n\} $, and assume the following:
	\begin{itemize}
		\item $ \GW(B_{\PP_n}(f_0,s_n))\geq r_n^2 $ 
		\item $\|f_{(n)}-f_0\|_{\PP_n} \leq \delta_n$,  and $ r_n \leq 2\delta_n$, 
		$s_n \leq \delta_n$.
		\item $\GW(B_{\PP_n}(f_{(n)},t)) \leq w_n't$~ for some ~$ t \leq \frac{r_n^4}{16\delta_n^2w_n'}$.
	\end{itemize}
	Then, when $Y_i = f_{(n)}(x_i) + \xi_i$, the following holds with high probability
	\begin{equation*}
	\|\hat{f}_n - f_{(n)}\|^2_{\PP_n} \geq c \cdot 
	\min\left\{\frac{r_n^8}{\delta_n^4w_n^2},\frac{r_n^4}{\delta_n^2}\right\}. 
	\end{equation*}
\end{lemma}
\begin{proof}[Proof of Lemma \ref{Lem:Twoptsfixed}]
	By Theorem \ref{chat}, it is enough to show that $t_{f_{(n)}}$ is greater than $\frac{r_n^4}{16\delta_n^2w_n}$. Namely, the  functional $H_{f_{(n)}}(\cdot)$ attains its unique maximum on a $t_{f_{(n)}} \geq \frac{r_n^4}{16\delta_n^2w_n}$. First, using the second assumption and the upper bound on $\delta_n$, we see that
	\[
	B_{\PP_n}(f_0, s_n) \subseteq B_{\PP_n}(f_{(n)}, \|f_{(n)}-f_0\|_{\PP_n}+ s_n) \subseteq B_{\PP_n}(f_{(n)}, 2\delta_n).
	\]
	By the convexity of $\F$, the function $t\mapsto \GW(B_{\PP_n}(f_{(n)},t))/t$ is nonincreasing. Therefore, by the first assumption and the above inclusion, for all $t \leq 2\delta_n$
	\[
	\GW(B_{\PP_n}(f_{(n)},t)) \geq \frac{r_n^2}{2\delta_n}t.
	\]
	Now, by the last equation and the second assumption, we know that for $t_1 = \frac{r_n^2}{2\delta_n} \leq 2\delta_n$
	\begin{align*}
	H_{f_{(n)}}(\frac{r_n^2}{2\delta_n}) &= \GW(B_{\PP_n}(f_{(n)},\frac{r_n^2}{2\delta_n})) - \frac{\left(\frac{r_n^2}{2\delta_n}\right)^2}{2} 
	\geq \frac{r_n^2}{2\delta_n} \cdot \left(\frac{r_n^2}{2\delta_n}\right) - \frac{\left(\frac{r_n^2}{2\delta_n}\right)^2}{2} \geq \frac{r_n^4}{8\delta_n^2}.
	\end{align*}
	
	Finally, by the last assumption, we know that for $t_2 = \frac{r_n^4}{16\delta_n^2w_n'} < t_1$ the following holds:
	\[
	H_{f_{(n)}}(\frac{r_n^4}{16\delta_n^2w_n'}) \leq \GW(B_{\PP_n}(f_{(n)}, \frac{r_n^4}{16\delta_n^2w_n'})) \leq w_n'\frac{r_n^4}{16\delta_n^2w_n'} = \frac{r_n^4}{16\delta_n^2}.
	\]
	Since $H_{f_{(n)}}(\cdot)$ is concave in $ t$, we conclude by 
	the last two equations that $ t_{f_{(n)}} \geq \frac{r_n^4}{16\delta_n^2w_n'}$ and the claim follows from Theorem \ref{chat}.
\end{proof}

In addition to Lemma \ref{Lem:Twoptsfixed}, we need the following two standard facts (which can be found, for example, in \cite{koltchinskii2011oracle}) to prove Theorem \ref{Thm:LowerBirgeMassartFixed}. 

\begin{lemma}[Sudakov Minoration]\label{Lem:sud}
	There exists a universal positive constant $C$ such that the following holds for any class $\F$ of real-valued functions:
	\begin{equation*}
	\GW(\F) \geq  \frac{C}{\sqrt{n}} \sup_{\epsilon > 0}
	\left\{\epsilon \sqrt{\log \Ncal(\epsilon,\F, \PP_n)} \right\}. 
	\end{equation*}
\end{lemma}
\begin{lemma}[Dudley Integral]\label{Lem:Dud}
	There exists a universal positive constant $C$ such that the following holds for any class $\F$ of real-valued functions:
	\begin{equation*}
	\GW(\F) \leq C\inf_{\epsilon > 0} \left(\epsilon + \frac{1}{\sqrt{n}}\int_{\epsilon}^{\mathrm{diam}(\F,\PP_n)/2} 
	\sqrt{\log \Ncal(u,\F, \PP_n)}du\right), 
	\end{equation*}
	where $\mathrm{diam}(\F,\PP_n)$ is the diameter of $\F$ with respect to $L^2(\PP_n)$.
\end{lemma}
\begin{proof}[Proof of Theorem \ref{Thm:LowerBirgeMassartFixed}]
	We aim to apply Lemma \ref{Lem:Twoptsfixed}. In order to satisfy its first assumption, we need to estimate the Gaussian width of $B_{\PP_n}(f_0, C_3n^{-1/p}\log(n)^{\gamma})$. By Lemma \ref{Lem:sud} and Eq. \eqref{A2} 
	\begin{equation}\label{Eq:51}
	\begin{aligned}
	\GW(B_{\PP_n}(f_0, C_3n^{-1/p}\log(n)^{\gamma})) &\geq \frac{C}{\sqrt{n}} \sup_{\epsilon > C_1n^{-1/p}\log(n)^{\gamma}}
	\left\{\epsilon \sqrt{\log \Ncal(\epsilon,B_{\PP_n}(f_0, C_3n^{-1/p}\log(n)^{\gamma},\PP_n)} \right\} \\&\geq  \frac{c_2}{\sqrt{n}}n^{-\frac{1}{p}}\log(n)^{\gamma}\sqrt{\log \Ncal(\epsilon,B_{\PP_n}(f_0, C_3n^{-1/p}\log(n)^{\gamma}),\PP_n)} \\&\geq c_3(p)n^{-\frac{1}{p}}\log(n)^{-\gamma(\frac{p}{2}-1)}.
	\end{aligned}
	\end{equation}
	Now, we aim to upper bound $\GW(B_{\PP_n}(f_{(n)},t))$ when
	$$   t = c_1(p)n^{-\frac{1}{2p}} \log (n)^{-\frac{\beta}{p} -\frac{3\gamma}{2}(\frac{p}{2}-1)},$$
	for some suitable $c_1(p)$. By using  Lemma \ref{Lem:Dud}, and Eq. \eqref{A3}, we see that
	\begin{equation}\label{Eq:blablabla}
	\begin{aligned}
	\GW(B_{\PP_n}(f_{(n)},t)) &\leq C\inf_{\epsilon >0} \left(\epsilon + \frac{1}{\sqrt{n}}\int_{\epsilon}^{t}\sqrt{\log\Ncal(u,B_{\PP_n}(f_{(n)}, t),\PP_n)}du\right)
	\\& \leq C\inf_{\epsilon >0 } \left(\epsilon + \frac {C_2(p)}{\sqrt{n}}\int_{\epsilon}^{t}\left(\frac{tn^{\frac 1{2p}}\log(n)^{\frac{\gamma}{2}(\frac{p}{2}-1)}\log(u^{-1})^{\frac{\beta}{p}}}{u}\right)^{\frac p2}du\right)
	\\& \leq C\inf_{\epsilon >0} \left(\epsilon + \frac {C_3(p)}{\sqrt{n}}\left(\frac{tn^{\frac 1{2p}}\log(\epsilon^{-1})^{\frac{\beta}{p}}\log(n)^{\frac{\gamma}{2}(\frac{p}{2}-1)}}{\epsilon}\right)^{\frac p2}\epsilon\right)
	\\&\leq C_4(p)n^{-\frac{1}{2p}}\log(n)^{\frac{\beta}{p} + \frac{\gamma}{2}(\frac{p}{2}-1)}t,
	\end{aligned}
	\end{equation}
	where we used the fact that $p > 2$ and set $$\epsilon =   C_1(p)n^{-\frac{1}{2p}}\log(n)^{\frac{\beta}{p} + \frac{\gamma}{2}(\frac{p}{2}-1)}t = C_2(p)n^{-\frac{1}{p}}\log(n)^{-\gamma(\frac{p}{2} -1 )} .$$ Hence, it is enough to assume \eqref{eq:weaker_cond} in place of \eqref{A41}.  Finally, we can apply Lemma \ref{Lem:Twoptsfixed} with the following parameters:
	$$ \|f_{(n)} - f_0\|_{\PP_n} \leq c_1(p)n^{-\frac{1}{2p}}\log(n)^{-\frac{\gamma}{2}(\frac{p}{2}-1)} =: \delta_n,$$
	where we used Eq. \eqref{A41}. Also, we set $s_n =  C_5(p)n^{-1/p}\log(n)^{\gamma}$,
	$$r_n^2 =  \min\{c_3(p),c_1(p)/8\}n^{-\frac{1}{p}}\log(n)^{-\gamma(\frac{p}{2}-1)}$$ and $$w_n =C_4(p) n^{-\frac{1}{2p}}\log(n)^{\frac{\beta}{p} + \frac{\gamma}{2}(\frac{p}{2}-1)}.$$ Therefore, Lemma \ref{Lem:Twoptsfixed} gives that
	\[
	\Ri_{\PP}(\hat{f}_n,\F) \geq c_5(p)\left(\frac{n^{-\frac{1}{p}}\log(n)^{-\gamma(\frac{p}{2}-1)}}{n^{-\frac{1}{2p}}\log^{\frac{\beta}{p} + \frac{\gamma}{2}(\frac{p}{2}-1)}(n)}\right)^2 \gtrsim n^{-\frac 1p}\log(n)^{-\frac{2\beta}{p} -3\gamma(\frac{p}{2}-1)},
	\]
	and the claim follows.
\end{proof}
\subsection{Proof of Corollary \ref{Thm:LowerBirgeMassartRandom}}\label{SectionThm}
For completeness, we state a consequence of \cite[Thm 5.11]{van2000empirical} (see Appendix for a proof). 
\begin{lemma}\label{Lem:BMbound}
	Let $\cG$ be a family of functions uniformly bounded by one, and $\PP$ be some probability measure. For any $f \in \cG$ and $t \geq 0$, let
	$\varepsilon(f,t)$ be the stationary point of
	$$ \frac{1}{\sqrt{n}}\int_{\frac{\delta}{c}}^{t}\sqrt{\log\Ncal_{[]}(u, B_{\PP}(f,t),\PP)}du = \delta$$
	for some absolute constant $c > 0 $. Then, for any $\varepsilon \geq \varepsilon(f,t) $ the following holds with probability of at least $1-C\exp(-cn\varepsilon^2)$ 
	
	\[
	\GW(B_{\PP}(f,t)) \lesssim \varepsilon + n^{-1/2}.
	\]
	
\end{lemma}

\begin{proof}[Proof of Corollary~\ref{Thm:LowerBirgeMassartRandom}]
	First, we take $m = C(p)\sqrt{n}\log(n)^{\frac{p\gamma}{2}(\frac{p}{2}-1)}$ in \eqref{A3} and set $f_{(n)}:=f_{(m)}$. By using Eq. \eqref{Eq:Distortion} and the estimate $d_n \lesssim(\log n)^{2} n^{-1/p}$  from \citep{rakhlin2017empirical}, we see that with high probability,
	\[
	\|f_{(n)} - f_0\|_{\PP_n} \leq 2\|f_{(n)} - f_0\|_{\PP} + C\log^2(n)n^{-\frac{1}{p}} \lesssim n^{-\frac{1}{2p}}\log(n)^{\frac{\gamma}{2}(\frac{p}{2}-1)}.
	\]
	Thus Eq. \eqref{A31} holds for $f_{(n)}$. 
	Next, by using Eq. \eqref{Eq:Distortion}, it is easy to see that Eq. \eqref{A21} hold with $\gamma = 2$. It remains to show that \eqref{eq:weaker_cond} holds. 
	
	First, by Eq. \eqref{Eq:Distortion}, we know that for $t \geq Cd_n$ the following holds:
	\[
	\GW(B_{\PP_n}(f_{(n)},t)) \leq \GW(B_{\PP}(f_{(n)},2t + d_n))  \leq \GW(B_{\PP}(f_{(n)},4t)) .
	\]
	To estimate the last term, we use our bracketing assumption along with the uniform boundedness of functions in $\Fcal$. These assumptions let us apply Lemma \ref{Lem:BMbound} that yields an upper bound on this term. In order to upper bound $\varepsilon(f_{(n)},4t)$ (the fixed point of Lemma \ref{Lem:BMbound}), we use Eq. \eqref{A4} and derive that
	\begin{align*}
	\frac{1}{\sqrt{n}}\int_{\frac{\varepsilon}{c}}^{4t} \left(\frac{tn^{\frac 1{2p}}\log(n)^{\frac{\gamma}{2}(\frac{p}{2}-1)}\log(u^{-1})^{\frac{\beta}{p}}}{u}\right)^{\frac p2}du 
	\lesssim \frac {1}{\sqrt{n}}\left(\frac{tn^{\frac 1{2p}}\log(n)^{\frac{\beta}{p}}\log(n)^{\frac{\gamma}{2}(\frac{p}{2}-1)}}{\varepsilon}\right)^{\frac p2}\varepsilon
	\end{align*} 
	where in the last inequality we used the fact that $p > 2$ and $\varepsilon\geq 1/n$. Setting the right-hand side equal to $\varepsilon$, we find an upper bound on the fixed point 
	$$ \varepsilon(f_{(n)},4t) \lesssim   n^{-\frac{1}{2p}}\log^{\frac{\beta}{p} + \frac{\gamma}{2}(\frac{p}{2}-1)}(n)t. $$  
	We now set $t \sim n^{-\frac{1}{2p}} \log (n)^{-\frac{\beta}{p} -\frac{3s_{p,\gamma}}{2}}$ and apply Lemma \ref{Lem:BMbound} with $\varepsilon \sim n^{-1/p}\log(n)^{-\gamma (p/2-1)} $, concluding that with probability of at least $1-C\exp(-c\Tilde{\Omega}(n^{1-2/p}))$ 
	\[
	\GW(B_{\PP}(f_{(n)},4t)) \lesssim \varepsilon \sim n^{-\frac{1}{2p}}\log^{\frac{\beta}{p} + \frac{\gamma}{2}(\frac{p}{2}-1)}(n)t.
	\]
	Hence, all the assumptions of Theorem \ref{Thm:LowerBirgeMassartFixed} hold with high probability. Therefore, by Eq. \eqref{Eq:Distortion} and Theorem \ref{Thm:LowerBirgeMassartFixed}, we conclude that 
	\begin{align*}
	\|\hat{f}_n - f^{*}\|_{\PP}^2 &\geq \frac{1}{4}\|\hat{f}_n - f^{*}\|^2_{\PP_n} - 2d_n^2 \geq \frac{1}{8}\|\hat{f}_n - f^{*}\|_{\PP_n}^2 \\&\gtrsim n^{-\frac 1p}\log(n)^{-\frac{2\beta}{p} -3\gamma(\frac{p}{2}-1)},
	\end{align*}
	and the claim follows.
\end{proof}
\subsection{Proof of Corollary \ref{Cor:SuppRandom}}
Throughout this subsection $\PP$ denotes the uniform measure on the sphere. In this proof we will use the following auxiliary lemmas:
\begin{lemma}[\cite{bronshtein1976varepsilon,dudley1999uniform}]\label{Lem:Bro}
	The following holds for all  $0 < \epsilon < 1 $, 
	\[
	2^{-d}\epsilon^{-\frac{d-1}{2}}\leq \log\mathcal{N}_{2}\left(\epsilon,\C,\PP \right) \leq \log \mathcal{N}_{\infty}\left(\epsilon,\C,\PP\right)
	\leq Cd^{5/2} \epsilon^{-\frac{d-1}{2}},
	\]
	where $\log \mathcal{N}_{\infty}\left(\epsilon,\C,\PP\right)$ denotes covering with respect to $\|\cdot\|_{\infty}$. 
\end{lemma}
\begin{lemma}[\cite{bronshtein1976varepsilon,dudley1999uniform}]\label{Lem:BroLower}
	Let $h_{B_d}$ be the support of function of the unit ball in $\R^d$, that is $h_{B_d} \equiv 1$. 
	Then for any $0 \leq \epsilon \leq c$, there exists a set of cardinality $ M(\epsilon,d) \geq 2^{10^{-d}\epsilon^{-(d-1)/2}}$, denoted by $h_1,\ldots,h_{K_{M(\epsilon,d)}}$, that has the following properties:
	\begin{itemize}
		\item $\forall $ $1 \leq i < j \leq M(\epsilon,d)$: $ \|h_{K_i} - h_{K_j}\|_{\PP} \geq \epsilon 	$.
		\item $\forall $  $1 \leq i \leq M(\epsilon,d)$: $\|h_{K_i}-h_{B_d}\|_{\PP} \leq 8\epsilon$.
	\end{itemize}
\end{lemma}
\begin{definition}
	We say that $f: \R^d \to \R$ is a \emph{$k$-piecewise simplicial linear} if $f$ is a convex piecewise linear function and its support can be written as a union of $k$-simplicials, and  $f$ is linear in each of them.        
\end{definition}
For the following lemma, let $\mathcal{P}_d$  be the regular simplex with $d+1$ faces that contains the unit Euclidean ball, with vertices at distance $d$ from the origin. 
\begin{lemma}[Theorem 4.4 in \cite{KurConv} ]\label{Lem:convEnt}
	Let $\F(\Gamma,\mathcal{P}_d)$ be the family of all the convex functions that are uniformly bounded by $\Gamma > 0$, and their domain is a simplex.
	Let $f_k \in \F(\Gamma,\mathcal{P}_d)$ be $k$ piecewise simplicial linear. Then, the following bound holds for $0 \leq \epsilon \leq t$:
	\[
	\mathcal{N}_{[]}(\epsilon, B_{\mathrm{Unif}}(f_k,t),\mathrm{Unif}) \leq C(d)\left(t/\epsilon\right)^{d/2}\log^{d+1}\left(\eps^{-1}\right)\log(\Gamma)
	\]
	where $B_{\mathrm{Unif}}(f_k,t) \subset \F(M,\mathcal{P}_d)$ and $\mathrm{Unif}$ denotes the Lebesgue volume measure.
\end{lemma}

In order to prove Corollary \ref{Cor:SuppRandom}, we first prove that when we restrict the LSE to $\mathcal{C}(C_1 \cdot d)$ (here $C_1$ is some absolute constant) we achieve the desired bound:
\begin{lemma}\label{fr}
	For every $d \geq 6, n \geq C(d)$, the following holds
	\begin{equation}\label{fr.eq}
	\argmin_{K \in \mathcal{C}(C_1 \cdot d)} \sum_{i=1}^n \left(Y_i - h_K(X_i) \right)^2 \geq c(d)\log^{\gamma_d}(n)\nnt,
	\end{equation}
	whenever $K^{*} \in \C$.
\end{lemma}
Then, in Subsection \ref{Sec:reduc}, we prove that with high probability,
\[
\argmin_{K \in \mathcal{C}} \sum_{i=1}^n \left(Y_i - h_K(X_i) \right)^2 = \argmin_{K \in \mathcal{C}(C_1 \cdot d )} \sum_{i=1}^n \left(Y_i - h_K(X_i) \right)^2
\]
whenever $K^{*} \in \C$, and therefore Corollary \ref{Cor:SuppRandom} follows.

\begin{proof}[Proof of Lemma \ref{fr}]
	Let $h_{K} \in \C $ and $t > 0$, denote by
	\[
	B_{\PP}(h_{K},t) := \{L \in \mathcal{C}(C_1 \cdot d): \|h_{K} - h_{L} \|_{\PP} \leq t \}.
	\]
	We will need the following Lemma (its proof appears in the next subsection).
	\begin{lemma}[The regular polytope lemma]\label{Lem:Regu}
		For every $m \geq  4^d,d \geq 2 $ there exists a polytope $P_{m,d} \subset B_d$ with $m$ vertices that satisfies the following: 
		\begin{enumerate}
			\item  $ \|\hp(x) - 1\|_{\infty} \leq C_1 m^{-\frac{2}{d-1}}$
			\item   $	\log \mathcal{N}_{[]}(\epsilon, B_{\Unif}( \hp,t),\Unif) \leq C(d)m\log(\epsilon^{-1})^{d} 
			\left(\frac{t}{\epsilon} \right)^{(d-1)/2} \ \ $
			for all $  0 \leq \epsilon \leq t $. 
		\end{enumerate}	
	\end{lemma} 
	Now, observe that the family is uniformly bounded by $C_1\cdot d$. Therefore, we can invoke Corollary \ref{Thm:LowerBirgeMassartRandom} with $f_{(m)} = \hp$, $\beta = d$, $f_0 = h_{B_d}$, and $p = \frac{d-1}{2}$. By Lemmas \ref{Lem:Bro}, \ref{Lem:BroLower}, and \ref{Lem:Regu}, we satisfy all the corollary's assumptions (observe that an upper bound on $\infty$-covering  implies that Koltchinski-Pollard entropy satisfies the same upper bound). Thus  $\Ri(\hat{h}_{\hat{K}_n},\C) \geq \Tilde{\Omega}(n^{-\frac{2}{d-1}}),$
	and the claim follows.
\end{proof}
\subsection{Proof of Lemma \ref{Lem:Regu}}
In this proof, we denote by $B_G(x,r)$ a geodesic ball on the sphere with center $x$ and radius $r$. Let $ \V:=\{v_i\}_{i=1}^{s} $ be a set of vertices of size $ m \leq s \leq 400d\ln(d)m$ with norm one. This set is a $m^{-\frac{1}{d-1}}$-net on $ \Sd $ with the following property: In any geodesic ball with radius  of $m^{-\frac{1}{d-1}}$ on the sphere there are at most $400d\ln(d)$ points. When $m \geq C^d$ such a set exists by \citep{boroczky2003covering}.

Now, we define a polytope $P_{m,d} := \mathrm{conv}\{\V\}$. First, observe that for each $x \in \R^d$ the following holds: 
\begin{equation}\label{Eq:support}
\hp(x) =\max_{v \in \V}\|x\|\|v\| \cos(\angle(x,v_{\pi}))= \|x\|\cos(\angle(x,v_{\pi})),
\end{equation}
where $v_{\pi}$ is the closest point in the net to $x/\|x\|$. Then, using the fact that $\cos(t) = 1 - t^2/2 +O(t^4),$ we know that when $m \geq C_1^d$, for all $u \in \Sd$ 
$$
0 \leq 1-\hp(u) \leq C\mmt
$$
and the first part of the Lemma follows. 

Now, we prove the last part of  Lemma \ref{Lem:Regu}, that is an upper bound on the entropy numbers. 
Let us consider the support function when it is restricted to the facets of regular simplex $\mathcal{P}_d$, denoted by $S_k$, $ 1 \leq k \leq d+1 $. Now we prove the following:
\begin{lemma}
	For any $k\in\{1,\ldots,d+1\}$,  $\hp:S_k \to \R$ is at most $C\ln(d)m$ piecewise linear, and every piece has at most $C(d\ln(d))^24^{d-1} \ $ $(d-2)$- facets.  
\end{lemma}
\begin{proof}
	\begin{figure}[htbp]
		\centering
		\includegraphics[width=0.5\linewidth]{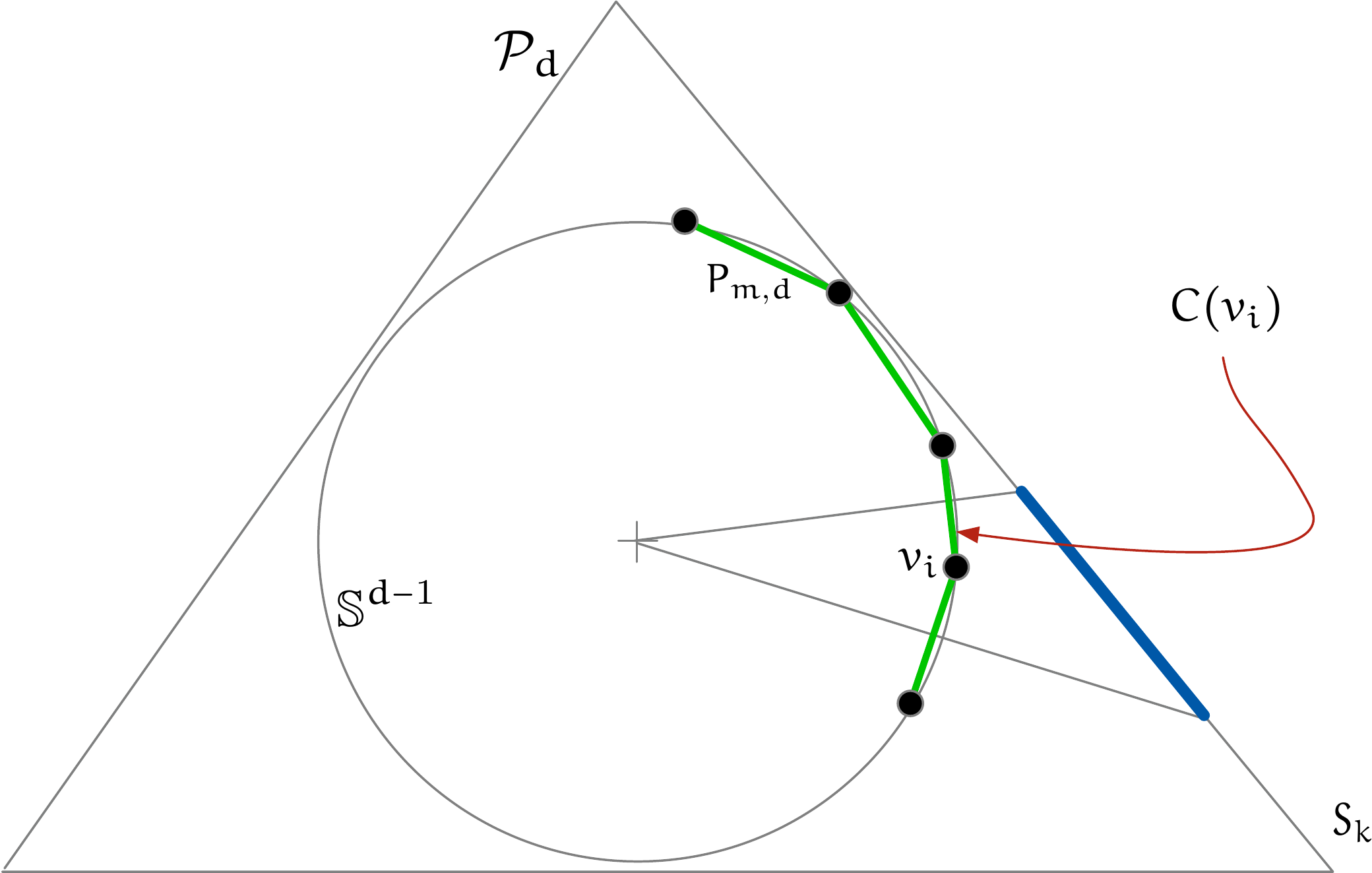}
		\caption{Illustration of the proof. The function $\hp$ is linear on the set $R^{-1}(C(v_i))$ (in blue) and convex piece-wise linear on $S_k$.}
		\label{fig}
	\end{figure}
	
	First, we use the Voronoi cells of the vertices $\V$ when they are restricted to $\Sd$, that is each cell, denoted by $C(v_i)$, $1 \leq i \leq |\V|$, is defined by $$ C(v_i) := \{x \in \Sd: \|v_i -x \| \leq \|v_j -x \| \ \ \forall   j \ne i\}.$$ 
	By Eq. \eqref{Eq:support} we know that for each $x \in \R^{d}$ the value of $\hp$ is attained on the closest $v_{i} \in \V$ to $x/\|x\|$ (since all the vertices of $ P_{n,d} $ have the same norm). Thus, we can write $$S_k = \bigcup_{i=1}^{|\V|}R^{-1}(C(v_i)) \cap S_k $$ 
	where $R$ denotes the radial function from $\mathcal{P}_d$ to $\Sd$. Moreover, observe that $$\hp: R^{-1}(C(v_i))(x) \equiv v_i^Tx.$$
	Since $\hp:S_k \to \R$ is a convex function, we conclude that it i also piecewise linear. Due to the regularity of the net  $\V$, we know that number of pieces can be bounded by $2|\V|/(d+1) \leq C\ln(d)m$. 
	
	Next, observe that the number of $d-2$ facets of each piece, which corresponds to some $v_i \in \mathcal{V}$, is determined by the number of neighbors of the Voronoi cell $C(v_i)$. By the construction of $ \V $, it can be bounded by $C(d\ln(d))^24^{d-1}$.
	To see this, since $\mathcal{V}$ is a $m^{-1/(d-1)}$-net, clearly all the neighbors of $v_i \in \mathcal{V}$ lie in $B_G(v_i,4\eps)$. Moreover, we can bound the number of vertices in this ball by 
	\begin{align*}
	|\{v \in \mathcal{V}: v \in B_G(v_i,4\eps)\}| &\leq 400d\ln(d) \cdot 
	|\Ncal(\eps,B_G(v_i,4\eps),\mathrm{Unif})| \\& \leq  C(d\ln(d))^24^{d-1}, 
	\end{align*}
	and the claim follows.

\end{proof}
Now, we are ready to prove Lemma \ref{Lem:Regu}. It is easy see that each $\hp:S_k \to \R$ is uniformly bounded by $d$ (since $K \subset B_d$ implies that $h_K$ restricted to $S_k$ is bounded by $d$). Moreover, since we assume that all $h_{L} \in B_{\PP}(h_K,t)$ are bounded in $C_1d$, we also know that our family is uniformly bounded by $Cd^2$. Now, recall that all of $\hp \ \ $   $C\sqrt{m}/(d-1)$-pieces have at most $(Cd\ln(d))^24^{(d-1)}$ facets, hence, it is also $C(d)\sqrt{m}$ piecewise simplicial linear.  Thus, we can apply Lemma \ref{Lem:convEnt} and  find a $(d+1)^{-1/2}\epsilon$-net, denoted by $\{f_{k,j}\}_{j=1}^{S_{\epsilon}}$,  with respect to $\mathrm{Unif}(S_k)$, where 
\[
S_{\epsilon}:= \log\mathcal{N}_{[]}((d+1)^{-1/2}\epsilon, B_{\mathrm{Unif}(S_k)}( \hp,t),\mathrm{Unif}(S_k)) \leq
C(d)m\log(\epsilon^{-1})^{d} \left(\frac{t}{\epsilon} \right)^{(d-1)/2}.
\] 
Recall that we aim to bound the entropy numbers of the the support function on the sphere. Therefore, we use the radial function to project each facet $S_k$ onto the sphere, and show that the set of functions  $$\{\|R^{-1}(x)\|^{-1}f_{k,j}(R^{-1}(x))\}_{j=1}^{S_{\epsilon}}$$ forms an $(d+1)^{-1/2}\epsilon$-net on $B_{\mathrm{Unif}(R(S_k))}(\hp,t)$ with respect to uniform measure.

Let $h_{K} \in  \mathcal{C}(C_1 \cdot d) $ (which is convex and bounded by one) that is at least $(d+1)^{-1/2}\epsilon$ far from $\hp$, when the support functions are restricted to $R(S_k)$. Denote by $h_{\pi,K}$  the closest member to the $h_{K}$ with respect to  $\mathrm{Unif}(S_k)$ in the aforementioned set. Then,
\begin{align*}
&\sqrt{\int_{R(S_k)}(h_K(x) - \|R^{-1}(x)\|^{-1}h_{\pi,K}(R^{-1}(x))^2dS(x)} 
\\&= \sqrt{\int_{R(S_k)}\|R^{-1}(x)\|^{-2}(h_K(R^{-1}(x)) - h_{\pi,K}(R^{-1}(x)))^2dS(x)} \\&=\sqrt{\int_{S_k}(h_K(x) - h_{\pi,K}(x))^2\frac{u_i^Tx}{\|x\|^{d+2}}dx} \leq \sqrt{\int_{S_k}(h_K(x) - h_{\pi,K}(x))^2dx} \leq (d+1)^{-1/2}\epsilon,
\end{align*}
where we used the homogeneity of the support function, and Lemma \ref{Lem:CV}. Therefore, we found our desired net.
Now, we can use  \citep[Lemma 2.8]{gao2017entropy}, and conclude that for all $0 \leq \epsilon \leq t$
\begin{align*}
\log \Ncal_{[]}(\epsilon, B_{\Unif}( \hp,t),\PP) &\leq   \sum_{k=1}^{d+1}\log \mathcal{N}_{[]}(\frac{\epsilon}{\sqrt{d+1}}, B_{\mathrm{Unif}(R(S_k))}( \hp,t),\mathrm{Unif}(R(S_k)))
\\&\leq   \sum_{k=1}^{d+1}\log \mathcal{N}_{[]}
(\frac{\epsilon}{\sqrt{d+1}}, B_{S_k}( \hp,t),\mathrm{Unif}(S_k))
\\&\leq C(d) m \log(\epsilon^{-1})^{d} \left(\frac{t}{\epsilon} \right)^{(d-1)/2},
\end{align*}
and the claim follows.
\subsection{Reduction to Lemma  \ref{fr}}\label{Sec:reduc}
We will show that when $K^{*} \in \C$, the LSE only considers convex sets in $\mathcal{C}(C_1 \cdot d )$. First, observe that the score of $h_{K^*}$ is bounded by $5$ when $n$ is large enough. To see this,
\begin{equation}\label{Eq:bla}
\begin{aligned}
n^{-1}\sum_{i=1}^{n}(Y_i - h_{K^{*}}(X_i))^2 &\leq 2n^{-1}\left(\sum_{i=1}^{n}\xi_i^2+ h_{K^{*}}(X_i)^2\right) \\&\leq 2\left(\int_{\Sd}h_K^2(x)d\PP(x) + 1 + O(\frac{1}{\sqrt{n}})\right) \leq 5.
\end{aligned}
\end{equation}
where we used the fact that $h_{K^*}(x) \leq 1$ when $x \in \Sd$.

Now, let $K \notin \mathcal{C}(C_1 \cdot d )$. Then, $K$ has a vertex (denoted by $v_1$)  with $\|v_1\|_2 > C_1d$. In this case for $C_1$ that is large enough, there exists a spherical cap of $\Sd$, denoted by $C(v_1)$, with center $v_1/\|v_1\|$ and normalized surface area of  $1/3$, such that $$h_{K}(u) \geq 100 \ \ \forall u \in C(v_1).$$ 

To see this, observe that if $C(v_1)$ has a normalized surface of $1/3$, then it implies that the geodesic distance between $v_1/\|v_1\|$ and its boundary is at most $\pi/2-c_2/(d-1)$ (for some fixed $c_2 \geq 0$). Thus, for any $u \in C(v_1)$ the following holds: 
\begin{align*}
h_{K}(u) &\geq \|v_1\|\cos(\angle(v_1,u)) \geq C_1d\cos(\pi/2-c_2/(d-1))  \\&= C_1d\sin(c_2/(d-1)) \geq C_1c_2/2,
\end{align*}
therefore if we choose $C_1$ to be large enough then $h_{K}(u) \geq 100$ for all $u \in C(v_1)$.

Now, using the fact that $X_i$ are i.i.d uniform on the sphere, we know that with  probability of at least $1-e^{-cn}$, there are at least $n/4$ points in this cap (using concentration for Bernoulli random variables). Moreover, since $\xi_i \sim N(0,1)$, we know that with high probability, at least $0.49$ of the points that lie in this cap, the $\xi_i$ are negative. Thus, we conclude that with high probability
\begin{align*}
n^{-1}\sum_{i=1}^{n}(Y_i - h_{K}(X_i))^2 &\geq n^{-1}\sum_{X_i \in C(v_1),~ \xi_i \leq 0}(h_{K}(X_i)- 1 -\xi_i)^2 \geq 11.
\end{align*}
Using the last equation and Eq. \eqref{Eq:bla}, only sets in $\mathcal{C}(C_1 \cdot d)$ will be considered, and the claim follows.

\bibliography{bib}
\appendix
\section{Auxiliary results and basic facts}
\begin{proof}[Proof of Lemma \ref{Lem:BMbound}]
	Under the conditions of the Lemma, \cite[Thm 5.11]{van2000empirical} implies that for the fixed point $\varepsilon$,
	\[
	\sup_{g \in B_{\PP}(f,t)}\big| \PP_n(g) - \PP(g)\big| \lesssim \varepsilon,
	\]
	with probability of at least $1-C\exp(-c(n(\varepsilon/t)^2)) \geq 1-C\exp(-cn\varepsilon^2)$ since $t\leq 1$ due to uniform boundedness. Now, using de-symetrization argument (e.g. \citep{koltchinskii2011oracle}), we know that
	\[
	\E\GW(B_{\PP}(f,t)) - \frac{1}{\sqrt{n}} \lesssim \E\sup_{g \in B_{\PP}(f,t)}\big| \PP_n(g) - \PP(g)\big| \lesssim  \varepsilon.
	\]
	Therefore, in expectation we have the desired bound. Since the class in uniformly bounded by one, we can apply McDiarmid’s inequality which gives 
	\[
	\Pr\left(\GW(B_{\PP}(f,t)) - C(n^{-1/2} + \varepsilon) \geq u\right) \leq \exp(-cnu^2)
	\]
	and by setting $ u = C(\varepsilon + n^{-1/2})$ the claim follows. 
\end{proof}
\begin{lemma}[Jacobian of the radial function, see for example \cite{schneider2008stochastic}]\label{Lem:CV}
	Let $H= \{x \in \R^d\ : \  x^Tu = h\}$ be a $d-1$ hyper-plane and any integrable $f:H \to \R$, the following holds for the radial function $R(x) = x/\|x\|$:
	$$\int_{H}f(x)\frac{u^Tx}{\|x\|^{d}}dx = \int_{R(H)}f(R^{-1}(x))dS(x).$$
\end{lemma}

\begin{lemma}[Basic facts on the support function, see for example \citep{artstein2015asymptotic}]\label{Lem:BasF}
	\quad
	\begin{itemize}
		\item  For every $K \in \C$, the function $h_K$ can be extended to the whole of $\R^d$ via $h_K(x) = \|x\|h_K(x/\|x\|)$ and this extension makes $h_K$ a convex and 1-Lipschitz function on $\R^d$. 
		\item For any $\lambda > 0$, $K \in \mathcal{C}$ the following holds for all $x\in\R^d$: $h_{\lambda K}(x) = \lambda h_{K}(x)$.
	\end{itemize}
	
\end{lemma}

\end{document}